\documentclass[12pt]{amsart}
\usepackage{amssymb, mathrsfs,amsmath, amsxtra}
\usepackage{latexsym}
\usepackage{amsfonts}
\usepackage{shadow}
\date{}

\author[D. Alpay]{Daniel Alpay}
\author[H. Attia]{Haim Attia}
\address{(DA, HA) Department of Mathematics \newline
Ben Gurion University of the Negev \newline P.O.B. 653, \newline
Be'er Sheva 84105, \newline ISRAEL} \email{dany@math.bgu.ac.il,
atyah@bgu.ac.il}
\author[D. Levanony]{David Levanony}
\address{(DL) Department of Electrical Engineering \newline
Ben Gurion University of the Negev \newline P.O.B. 653, \newline
Be'er Sheva 84105, \newline ISRAEL }
\thanks{D. Alpay thanks the
Earl Katz family for endowing the chair which supported his
research. The research of the authors was supported in part by the
Israel Science Foundation grant 1023/07}
\email{levanony@ee.bgu.ac.il}

\title[On a class of Gaussian processes]
{On the characteristics of a class of Gaussian processes within
the white noise space setting} 
\begin{document}
\maketitle
\parindent 0cm
\newtheorem{Pa}{Paper}[section]
\newtheorem{Tm}[Pa]{{\bf Theorem}}
\newtheorem{La}[Pa]{{\bf Lemma}}
\newtheorem{Cy}[Pa]{{\bf Corollary}}
\newtheorem{Rk}[Pa]{{\bf Remark}}
\newtheorem{Pn}[Pa]{{\bf Proposition}}
\newtheorem{Pb}[Pa]{{\bf Problem}}
\newtheorem{Dn}[Pa]{{\bf Definition}}
\newtheorem{Ex}[Pa]{{\bf Example}}
\numberwithin{equation}{section}
\def\L{\mathbf L}
\def\R{\mathbb R}
\def\N{\mathbb N}
\def\C{\mathbb C}
\def\s{\mathscr S(\R)}
\def\ss{\mathscr S'(\R)}
\def\(s){\mathscr S(\R^2)}
\def\F{\mathcal F}
\def\P{\mathcal P}
\def\W{\mathcal W}
\def\Dom{{\rm dom}~(T_m)}
\def\Doms{{\rm dom}~(T_m^*)}
\def\Def{\stackrel{{\rm def.}}{=}}
\begin{abstract}
Using the white noise space framework, we define a class of
stochastic processes which include as a particular case the
fractional Brownian motion and its derivative. The covariance
functions of these processes are of a special form, studied by
Schoenberg, von Neumann and Krein.
\end{abstract}
\maketitle \tableofcontents

\section{Introduction}
Using the white noise framework, we define and
study Gaussian processes whose covariance functions are the form
\begin{equation}
K_r(t,s)=r(t)+r(s)^*-r(t-s)-r(0), \quad
t,s\in{\mathbb R},
\label{kernel}
\end{equation}
where $r$ is a continuous function such that
\begin{equation*}
r(-t)=r(t)^*, \quad
t\in{\mathbb R},
\end{equation*}
and the kernel $K_r(t,s)$ is positive (in the sense of
reproducing kernels) on the real line. As we will recall in the
sequel, such functions $r$ have been investigated for a long
time. Still, their applications in stochastic calculus seem to
have been only partially developed. We mention in particular the
recent work \cite[p. 103]{MR2338856}. In that work the notion of
processes with covariance measure is introduced, and stochastic
processes with covariance function of the form $K_r$ are shown to
belong to this class. We note however that the methods of
\cite{MR2338856} and of the present
paper are completely different.\\

Functions $r$ such that the kernel $K_r(t,s)$ is
positive have been characterized by von Neumann
and Schoenberg when $r$ is real; see
\cite[Theorem 1, p. 229]{MR0004644}. For the
case of complex--valued functions, see Krein
\cite{MR0012176} and \cite{MR0012175}, and
Akhiezer, \cite[pp. 267--269]{akhiezer_russian}
and the references therein. These are functions
of the form
\begin{equation}
\label{kipour_2007}
r(t)=r_0 +i\gamma t-\int_{\R}
\Big\{e^{itu}-1-\frac{itu}{u^2+1}
\Big\}\frac{d\sigma(u)}{u^2},
\end{equation}
where $r_0=r(0)$ and $\gamma$ are real numbers
with $d\sigma$ a positive measure on
${\mathbb R}$ defined by an increasing right continuous function
$\sigma$,
such that
\begin{equation}
\label{sigma!}
\int_{\mathbb
R}\frac{d\sigma(u)}{u^2+1}<\infty.
\end{equation}
That the form \eqref{kipour_2007} is sufficient to insure the
positivity of the kernel $K_r(t,s)$
follows from the formula
\begin{equation}
\label{qaz}
K_r(t,s)= \int_{\mathbb
R}\frac{e^{itu}-1}{u}\frac{e^{-isu}-1}
{u}d\sigma(u),
\end{equation}
see for instance \cite[Theorem 4 p.
115]{MR1321817}. For a discussion of this formula, see also
\cite[(7) p. 25]{lifs}. The idea of the proof of the
converse is given in the next section.\\

Note that
\[
K_r(t,s)=K_{r-r(0)}(t,s),
\]
and therefore one can always assume that $r(0)=0$.\\

In the real valued case, with $r(0)=0$, the function $r$ takes
the form
\begin{equation*}
r(t)=\int_{\mathbb R} \dfrac{1-\cos(tu)}{u^2}d\sigma(u).
\end{equation*}

The fractional Brownian motion then corresponds to the
choice
\[
d\sigma(u)=\frac{1}{2\pi}|u|^{1-2H}du,\quad H\in(0,1),
\]
giving
\begin{equation}
\label{mr}
r(t)=\frac{V_H}{2}|t|^{2H}, \quad {\rm with}\quad V_H=
\frac{\Gamma(2-2H)\cos(\pi H)}{\pi(1-2H)H},
\end{equation}
and
$$
K_r(t,s){\stackrel{{\rm def.}}{=}}
k_H(t,s)=\frac{V_H}{2}(|t|^{2H}+|s|^{2H}-|t-s|^{2H}),$$
with $\Gamma$ denoting the Gamma function, as can be seen using
the formulas
\[
\begin{split}
\int_{\mathbb R}\frac{1-\cos(tu)}{u^{2H+1}}&=-2|t|^{2H}\cos(\pi
  H)\Gamma(-2H)\\
\int_{\mathbb R}\frac{1-\cos(tu)}{u^{2}}&=\pi |t|.
\end{split}
\]
When furthermore $H=1/2$, then $V_H=1$, $r(t)=\left|t\right|/2$
and for $t,s\ge 0$, $K_r(t,s)={\rm min}(t,s)$.\\

By a theorem of Kolmogorov there exists a Gaussian stochastic
process $\{B_H(t)\}$ indexed by ${\mathbb R}$, which is called the
fractional Brownian motion (with Hurst parameter $H\in(0,1)$),
such that
$$k_H(t,s)=E(B_H(t)B_H(s)),\quad   t,s\in\R.$$
Stochastic calculus for $B_H$ has been developed for quite some
time; see for instance \cite{DS}, \cite{MR1741154}, \cite{DVZ},
\cite{fbm88}. In a subsequent paper we show how most of the
results of these works are extended to the case of general
covariance functions of the form \eqref{kernel} above.\\

We mention also that
functions $r$ of the form \eqref{kipour_2007}
appear first in the work of Paul L\'evy \cite{MR1556734},
in the result characterizing characteristic
functions of infinitely divisible laws. More precisely,
the characteristic function of a random variable is infinitely
divisible if and only if, it is of the form $\exp r(t)$, where
$r(t)$ is of form \eqref{kipour_2007}; see
\cite[Representation Theorem]{MR0203748}.
Similarly,
$(Z_u)_{u\ge 0}$ is an infinitely divisible
random process if and only if
\[
E \left(e^{itZ_u}\right)=e^{-ur(t)},
\]
where $r$ is of the form \eqref{kipour_2007},
with $r(0)=0, \gamma=0$, and is called the {\sl
characteristic exponent of the L\'evy process}.
See \cite[p. 12]{bertoin}, \cite[formula (9), p.
353]{MR1556734}. This is the L\'evy-Khintchine
formula, see \cite[(formula (9) p.
353]{MR1556734}.
We also mention that positive kernels of
 the form $K_r$ appear in
the theory of Dirichlet spaces; see
\cite[p. 5-12]{deny}. These aspects of the theory of
kernels of the form $K_r$ will not be studied in the
present paper.\\

The paper consists of 7 sections including the introduction, and
its outline is as follows. In Section 2 we study and characterize
the reproducing kernel Hilbert space associated with a kernel
$K_r(t,s)$. In Section 3 we associate with certain kernels
$K_r(t,s)$ an operator which will play a key role in the
construction of the stochastic process with covariance
$K_r(t,s)$. This paper uses Hida's white noise space theory, and
in Section \ref{S4} we review the main features from white noise
space theory which we will subsequently use.
 In Section \ref{S5} we recall the definition
of the Wick product and of the Kondratiev space. Stochastic
processes with covariance function $K_r(t,s)$ are built in Section
6,
and their derivatives are studied in Section 7.\\


Some of the results presented here have been
announced in the note \cite{MR2414165}. The results of the present work are
to be used in a subsequent paper where stochastic analysis of the
processes considered here is developed.\\

Finally a word on notation. We denote the
Fourier transform by
\begin{equation}
\stackrel{\wedge}{f}(u)=\int_{\R}e^{-iux}f(x)dx.
\label{fourier}
\end{equation}
%
%
The inverse transform is then given by
\begin{equation}
\stackrel{\vee}{f}(u)=\frac{1}{2\pi}\int_{\R}e^{iux}f(x)dx.
\end{equation}

The same notations are used for the Fourier
transform and inverse Fourier transform of distributions.\\

We set
\[
{\mathbb N}=\left\{1,2,3,\ldots\right\}
\quad{\rm and}\quad
{\mathbb N}_0={\mathbb N}\cup\left\{0\right\},
\]
and denote by $\ell$ the set of sequences
\begin{equation}
\label{eq:ell}
(\alpha_1,\alpha_2,\ldots),
\end{equation}
indexed by ${\mathbb N}$  with values in
${\mathbb N}_0$, for which only a finite
number of elements $\alpha_j\not=0$.

\section{Some remarks on the kernels $K_r(t,s)$}
As already mentioned, real-valued functions $r$
for which the kernel $K_r(t,s)$ is positive on
the real line (in the sense of reproducing
kernels) were characterized by Schoenberg and von
Neumann. For complex-valued functions, and by
different methods, the following  theorem has been
given by Krein in 1944; see
\cite[Theorem 2, p.
256]{MR1321817}.
See also \cite[\S 9 p.
268]{akhiezer_russian}. In Krein's result, the
case $a=\infty$ is allowed, and then,
$t,s\in{\mathbb R}$.

\begin{Tm} \cite[Theorem 2, p.
256]{MR1321817}
The kernel $K_r(t,s)=r(t)+r(s)^*-r(t-s)-r(0)$ is
positive for $t,s\in[-a,a]$ if and only if it is
of the form \eqref{kipour_2007}:
\begin{equation*}
r(t)=r(0) +i\gamma t-\int_{\R}
\Big\{e^{itu}-1-\frac{itu}{u^2+1}\Big\}
\frac{d\sigma(u)}{u^2},\quad t\in[-a,a].
\end{equation*}
\end{Tm}
For completeness, let us recall that Akhiezer's
proof \cite[pp. 268-270]{akhiezer_russian} goes
along the following lines; one first shows that
$r$ satisfies an inequality of the form
\begin{equation}
|r(t)|\le M(1+|t|^3)
\label{ineq45}
\end{equation}
for some positive number $M$ (we recall the proof of this
inequality in the sequel; see Lemma \ref{111} below). One then
shows that the function $H(z)=z^2\int_0^\infty r(t)^*e^{itz}dt$,
which, in view of \eqref{ineq45}, is analytic in the open upper
half-plane, satisfies (see \cite[(2) pp. 268]{akhiezer_russian}
and \cite[p. 227]{kl2})
\[
\frac{H(z)-H(w)^*}{z-w^*}
=zw^*\iint_{\R_+^2}K_r(t,s)e^{itz}
e^{-isw^*}dtds,
\]
and in particular has a positive imaginary part in the open upper
half-plane.
To conclude the proof, one uses Herglotz's representation formula for
analytic functions with a positive imaginary part in the open
upper half-plane (see for instance \cite[Theorem 4.7, p.
25]{MR48:904}, \cite[Theorem 2 p. 220]{ak-vol1}).\\

Formula \eqref{qaz} allows to characterize the reproducing kernel
Hilbert space associated with $K_r$ in terms of de Branges
spaces. See Theorem \ref{tm:paris} below. We first make the
following remarks: Let $\chi_t$ denote the function of a real
variable $u$
\[
\chi_t(u)=\frac{e^{itu}-1}{iu},\quad t\in{\mathbb R},
\]
and let ${\mathscr M}_T$ denote the closed
linear span in ${\mathbf L}_2(d\sigma)$ of the
$\chi_t$ for $|t|\le T$. Assume that ${\mathscr
M}_T\not ={\mathbf L}_2(d\sigma)$. Then,
${\mathscr M}_T$ is a reproducing kernel Hilbert
space with reproducing kernel of the form
\begin{equation*}
\frac{A(T,\lambda)A(T,\omega)^{*}
-B(T,\lambda)B(T,\omega)^{*}}{-i(\lambda-\omega^*)}
\end{equation*}
where $A(T,\lambda), B(T,\lambda)$ are entire
functions of finite exponential type. See \cite{dbbook}, \cite{Dmk}.\\

\begin{Rk}
The spaces ${\mathscr M}_T$ were introduced by
de Branges and play a key role in prediction
theory. See \cite{dbbook, Dmk, dym-70}. When
$d\sigma(u) =du$ we have
$$
{\mathscr M}_T={\mathbf H}_2\ominus e_T {\mathbf H}_2,$$ where
${\mathbf H}_2$ denotes the Hardy space of the open half-plane,
and where $e_T(z)=e^{izT},B(T,\cdot)=1$ and $A(T,\cdot)=\chi_T$,
with the reproducing kernel given by
$$\frac{1-e_T(\lambda)e_T(\omega^*)^*}{-i(\lambda-\omega^*)}.$$
\end{Rk}

Let $\s$ denote the Schwartz space of rapidly decreasing
functions, and let $\ss$ denote the the topological
dual of $\s$, that is,  the space of tempered
distributions.
In view of the next two results, we recall the
following: Condition \eqref{sigma!} insures that
the measure $d\sigma$ has a Fourier transform
$\widehat{d\sigma}$ which is a tempered distribution.
Furthermore, this Fourier
transform induces a distribution
$\widehat{d\sigma}(t-s)$ on the Schwartz space $\(s)$ of functions of
two real variables via the formula
\begin{equation}
\label{f_t} \langle
\widehat{d\sigma}(t-s),\phi(t,s)\rangle=
\int_{\mathbb R}\left(\iint_{{\mathbb
R}^2}e^{-i(t-s)u}\phi(t,s)dtds\right)d\sigma(u).
\end{equation}

When $\int_{\mathbb R}d\sigma(u)<\infty$, we
have that
\[
\begin{split}
\langle\widehat{d\sigma}(t-s),\phi(t,s)\rangle&
=\iint_{\R^2}\left\{\int_\R
e^{-i(t-s)u}d\sigma(u)
\right\}\phi(t,s)dtds\\
&= \iint_{{ \mathbb
R}^2}\widehat{d\sigma}(t-s)\phi(t,s)dtds.
\end{split}
\]

\begin{Tm}
\label{distribution} Let $T<\infty$. The reproducing kernel
Hilbert space ${\mathcal H}_T(K_r)$ associated with $K_r(t,s)$ for
$t,s\in[-T,T]$ consists of functions of the form
\begin{equation}
\label{Ft}
F(t)=\int_{\mathbb R}
\frac{e^{itu}-1}{iu}f(u)d\sigma(u),\quad t\in[-T,T],
\quad f\in{\mathscr M}_T
\end{equation}
with norm
\begin{equation}
\label{Ft1}
\|F\|_{{\mathcal H}_T(K_r)}=
\|f\|_{{\mathbf L}_2(d\sigma)}.
\end{equation}
Moreover, by extending $F(t)$ to the real line by
formula \eqref{Ft}, $F$ defines a tempered
distribution and
\[
F^\prime(t)=2\pi\stackrel{\vee}{fd\sigma}(t)
\]
in the sense of distributions, and where
$\stackrel{\vee}{fd\sigma}$ denotes the inverse Fourier transform of the
tempered distribution defined by $fd\sigma$.
\label{tm:paris}
\end{Tm}
{\bf Proof:} Let $F$ be of the form \eqref{Ft}
and assume
$F\equiv 0$. Then $f$ is orthogonal to all
$\chi_t$, $t\in[-T,T]$ and
therefore $f=0$ since $f\in {\mathscr M}_T$.
Thus, \eqref{Ft1}
defines a norm. Furthermore, the choice
$f=(\chi_s)^*$ in \eqref{Ft} leads to
\[
\begin{split}
\left\langle F,K_r(\cdot, s)
\right\rangle_{{\mathcal H}_T(K_r)}
&=\left\langle f,(\chi_s)^*
\right\rangle_{\L_2(d\sigma)}\\
&=\int_{\mathbb R} f(u)\chi_s(u) d\sigma(u)=F(s),
\end{split}
\]
which proves the first claim.\\

We note that the function $F(t)$ extends to a
continuous function on ${\mathbb R}$. Indeed,
for every $t,h\in\R$ we have
\begin{equation}
\label{wsx}
\left|\frac{e^{i(t+h)u}-e^{itu}}{iu}\right|=
\left|\frac{e^{ihu}-1}{iu}\right|
\leq\begin{cases}
\left|h\right|\quad{\rm if}~~~\left|u\right|\leq 1,\\
~~\\
\dfrac{2}{\left|u\right|} \quad{\rm
if}~~~\left|u\right|>1.
\end{cases}
\end{equation}
Using these inequalities we use the dominated
convergence theorem to prove that
\[
\lim_{h\rightarrow 0}F(t+h)=F(t),\quad t\in{\mathbb R}.
\]
Furthermore, \eqref{wsx} leads to the bound:
\[
\begin{split}
|F(t)|&\le|t|\int_{-1}^1|f(u)|d\sigma(u)
+2\int_{|u|\ge
1}\left|\frac{f(u)}{u}\right|d\sigma(u).
\end{split}
\]
We recall that $\sigma$ is assumed right continous. When it
has a jump at $0$, we define
\[
F(t)=tf(0)(\sigma (0)-\sigma (0_-))+\int_{\mathbb R} f(u)\chi_s(u)
d\sigma_1(u),
\]
where $d\sigma_1$ has no jump at $0$.
The function $F$ is in particular slowly
growing, and therefore defines a tempered
distribution (see \cite[Th\'eor\`eme VI p.
239]{Schwartz66}, \cite[\S4, p. 110]{Barros}).\\

Let $\varphi\in\s$. The integral
\[
\int_\R f(u)\varphi(u)d\sigma(u) =\int_{\mathbb
  R}\frac{f(u)}{u+i}\left( (u+i)
\varphi(u)\right)d\sigma(u)
\]
exists, since $(u+i)\varphi(u)$ is bounded and
since $1/(u+i)\in {{\mathbf L}_2(d\sigma)}$.
Thus,
\[
\begin{split}
\int_\R
\left(\int_\R
\left|\frac{e^{itu}-1}{iu}f(u)
\varphi^{'}(t)\right|d\sigma(u)
\right)dt
&\le
\int_{\left|u\right|\leq 1} \left|f(u)\right|
d\sigma (u)\int_\R\left|t\varphi ^{'}(t)\right|dt+\\
&\hspace{-8mm} +2\int_{\left|u\right|> 1}
\left|\frac{f(u)}{u}\right| d\sigma (u)
\int_\R\left|\varphi ^{'}(t)\right|dt<\infty.
\end{split}
\]

Using Fubini's theorem, we have
\[ \int_\R\left\{\int_\R\frac{e^{itu}-1}{iu}f(u)
d\sigma(u)\right\}\varphi^{'}(t)dt
=\int_\R\left\{\int_\R\frac{e^{itu}-1}{iu}
\varphi^{'}(t)dt\right\}
f(u)d\sigma(u),
\]
and by integration by parts, we obtain
\[
\begin{split}
\int_\R\left\{\int_\R\frac{e^{itu}-1}{iu}
\varphi^{'}(t)dt\right\}f(u)d\sigma(u)
&=-\int_\R\left\{\int_\R e^{itu}
\varphi(t)dt\right\}f(u)d\sigma(u)\\
&=-2\pi \int_\R\stackrel{\vee}{\varphi}(u)f(u)d\sigma(u).
\end{split}
\]
Thus,
\[
\int_\R
F(t)\varphi^{'}(t)dt= -2\pi
\int_\R\stackrel{\vee}{\varphi}(u)f(u)d\sigma(u),
\]
and therefore we obtain on the one hand
\[
\left\langle F,\varphi^\prime\right\rangle=
-2\pi\left\langle
fd\sigma,\stackrel{\vee}{\varphi}\right\rangle=-2\pi\left\langle
\stackrel{\vee}{fd\sigma},\varphi\right\rangle.
\]
On the other hand,
\[
\begin{split}
\left\langle F,\varphi^\prime\right\rangle=
-\left\langle F',\varphi\right\rangle
\end{split}
\]
Thus, $F^\prime=2\pi \stackrel{\vee}{fd\sigma}$.
\mbox{}\qed\mbox{}\\

In preparation for the proof of Theorem \ref{1} below we now prove
inequality \eqref{ineq45}.

\begin{La}
\label{111}
Assume the kernel $K_r(t,s)$ to be positive
in ${\mathbb R}$. Then \eqref{ineq45} is in
force, that is
\begin{equation*}
|r(t)|\le M(1+|t|^3)
\end{equation*}
for some positive number $M$.
\end{La}

{\bf Proof:} We follow the arguments in \cite[pp.
264-265]{akhiezer_russian}, with slight
modifications. We first note that we may assume
that $r(0)=0$. The positivity of the kernel
$K_r(t,s)$ implies that the matrix
\[
\begin{pmatrix}
K_r(t,t)&K_r(t,-t)\\K_r(-t,t)&K_r(-t,-t)
\end{pmatrix}
\]
has a non-negative determinant. Therefore,
\[
|2r(t)-r(2t)|\le |2{\rm Re}~r(t)|,
\]
and thus
\begin{equation}
\label{ineq_r}
|r(2t)|\le 4|r(t)|.
\end{equation}
Let
\[
R(t)=\frac{|r(t)|}{1+\left|t\right|^3}.
\]
Then, \eqref{ineq_r} implies that
\begin{equation}
R(2t)\le\frac{4(1+\left|t\right|^3)}{1+8\left|t\right|^3}R(t).
\label{chemin_vert}
\end{equation}
Let $T_0\in{\mathbb R}_+$ be such that:
\[
|t|\ge
T_0\,\,\Longrightarrow\,\,
\frac{4(1+\left|t\right|^3)}{1+8\left|t\right|^3}\le 1.
\]
It follows from \eqref{chemin_vert} that $R(t)$ is bounded in
${\mathbb R}$ by an expression of the form $M(1+|t|^3)$ for some
$M>0$. In fact, since $r$ is continuous, one may take
\[
M=\max_{t\in[0,T_0]}|r(t)|.
\]

\mbox{}\qed\mbox{}\\

\begin{Tm}
It holds that
\[
\frac{\partial^2}{\partial t\partial s}K_r(t,s)
=r''(t-s)= \widehat{d\sigma}(s-t)
\]
in the sense of distributions. Furthermore,  for $\varphi\in\s$,
$\stackrel{\vee}{d\sigma}\ast\varphi$ is a function and it holds that
$\varphi\in\s$
\begin{equation}
\label{qw}
\left\langle d\sigma\widehat{\varphi},
(\widehat{\varphi})^*\right\rangle_{\ss,\s}
=\int_{\mathbb
  R}{\varphi(-u)^*}
  (\stackrel{\vee}{d\sigma}\ast\varphi)(u)du\\
\end{equation}
\label{1}
\end{Tm}

Before the proof, we make the following observation.
We note that the right hand side of \eqref{qw} can
formally be rewritten as
\[
\int_{\mathbb
  R}{\varphi(-u)^*}(\int_{\mathbb
R}\stackrel{\vee}{d\sigma}(u-v)\varphi(v)dv)du
\]
where in general
\[
\int_{\mathbb
R}\stackrel{\vee}{d\sigma}(u-v)\varphi(v)dv
\]
is not a real integral, but an abuse of notation.\\

In the proof of the theorem, use is made of
properties of the space ${\mathcal O}_M$ of
multiplication operator in $\s$ (also called
$C^\infty$ functions slowly decreasing at
infinity), see \cite[p. 275]{Treves67}, and of
the space ${\mathcal O}^\prime_C$ of
distributions rapidly decreasing at infinity,
see \cite[p. 315]{Treves67}. Recall
(\cite[Theorem 30.3, p. 318]{Treves67}) that the
Fourier transform is one-to-one from ${\mathcal
O}_M$ onto ${\mathcal O}_C^\prime$ and from
${\mathcal O}_C^\prime$ onto ${\mathcal O}_M$.\\

{\bf Proof of Theorem \ref{1}:}
Since $d\sigma$
defines a distribution in $\ss$, we have that
$\stackrel{\vee}{d\sigma}\in\ss$; see for instance
\cite[Theorem 25.6, p. 276]{Treves67}. Let
$\varphi\in\s$.  The convolution
$\stackrel{\vee}{d\sigma}\ast\varphi$ is a function, and
belongs to ${\mathcal O}_M$; see
\cite[p.248]{Schwartz66}. So, by \cite[Theorem
30.3, p. 318]{Treves67}
\begin{equation}
\widehat{(\stackrel{\vee}{d\sigma}\ast\varphi)}\in
{\mathcal O}^\prime_C.
\label{conv}
\end{equation}
Now we compute \eqref{conv} using \cite[Theorem
30.4 p. 319]{Treves67} with (in the notation of
that book) $S=\stackrel{\vee}{d\sigma}\in \ss$ and
$T=\varphi\in{\s}\subset{\mathcal O}^\prime_C$
(see \cite[Example 30.1, p. 315]{Treves67} for the
latter inclusion), and obtain
\[
\widehat{(\stackrel{\vee}{d\sigma}\ast\varphi)}=
d\sigma\widehat{\varphi},
\]
which is a measure. Using the fact that
${\mathcal O}^\prime_C\subset\ss$ (see \cite[p.
318]{Treves67}, \cite{Schwartz66}), we have
\begin{equation}
\left\langle d\sigma
\widehat{\varphi},\widehat{\varphi}^*\right
\rangle_{\ss,\s}=
\int_{\R}\left|\widehat{\varphi}\right|^2d\sigma.
\label{act}
\end{equation}
But for $\psi\in\ss$ and $\varphi\in\s$ we have:
\begin{equation*}
\left\langle\widehat{\psi},\varphi\right
\rangle_{\ss,\s}=
\left\langle\psi,\widehat{\varphi}\right
\rangle_{\ss,\s}.
\end{equation*}
Let $\phi\in\s$. Then the function
\[
\psi\quad:\quad u\mapsto\varphi(-u)^*
\]
is also in $\s$
and we have:
\[
\begin{split}
\int_{\mathbb
  R}\phi(-u)^*(\stackrel{\vee}{d\sigma}\ast\varphi)(u)du&=
\left\langle\stackrel{\vee}{d\sigma}\ast\varphi,\psi
\right\rangle_{\ss,\s}\\
&=\left\langle d\sigma\widehat{\varphi},
\widehat{\psi}\right\rangle_{\ss,\s}\\
&=
\left\langle d\sigma\widehat{\varphi},
\widehat{\varphi}^*\right\rangle_{\ss,\s},
\end{split}
\]
since $\widehat{\psi}=(\widehat{\varphi})^*.$
Thus, \eqref{act} can be rewritten as:
\begin{equation*}
\begin{split}
\left\langle d\sigma\widehat{\varphi},
\widehat{\varphi}^*\right\rangle_{\ss,\s}
=\int_{\mathbb
  R}\varphi(-u)^*(\stackrel{\vee}{d\sigma}\ast\varphi)(u)du.
\end{split}
\end{equation*}
Now let $\phi\in\(s)$: We have
\[
\begin{split}
\left\langle\frac{\partial^2}{\partial t\partial s}
K_r,\phi\right\rangle&=
\left\langle K_r,
\frac{\partial^2}{\partial t\partial s}
\phi\right\rangle\\
&=\iint_{\R^2}K_r(t,s)
\frac{\partial^2}{\partial t\partial s}
\phi(t,s)dtds\\
&=\iint_{\R^2}r(t)
\frac{\partial^2}{\partial t\partial s}
\phi(t,s)dtds\\
&\hspace{5mm}+\iint_{\R^2}r(s)^*
\frac{\partial^2}{\partial t\partial s}
\phi(t,s)dtds\\
&\hspace{10mm}-\iint_{\R^2}r(t-s)
\frac{\partial^2}{\partial t\partial
s}\phi(t,s)dtds.
\end{split}
\]
Since $r$ satisfies inequality \eqref{ineq45}, the above
integrals make sense; the first and the second integrals on the
right handside are identically zero since $\phi$ is a Schwartz
function. Thus, since
\[
\begin{split}
\left\langle
r''(t-s),\phi(t,s)\right\rangle&=-\left\langle
r(t-s),\frac{\partial^2}{\partial t\partial
s}\phi\right\rangle\\
&=-\iint_{\R^2}r(t-s)\frac{\partial^2}{\partial
t\partial s}\phi(t,s)dtds,
\end{split}
\]
it follows that
\[
\left\langle
\frac{\partial^2}{\partial t\partial s}
K_r,\phi\right\rangle
=\left\langle r''(t-s),\phi\right\rangle.
\]
Using \eqref{kipour_2007}
\[r(t-s)=
i\gamma (t-s)-\int_{\R}
\Big\{e^{i(t-s)u}-1-
\frac{i(t-s)u}{u^2+1}\Big\}\frac{d\sigma(u)}{u^2},
\]
we get
\[
\begin{split}
\iint_{{\mathbb R}^2}
r(t-s)\frac{\partial^2}{\partial t\partial s}
\phi(t,s)dtds &=\iint_{\R^2} i\gamma(t-s)
\frac{\partial^2}{\partial t\partial s}
\phi(t,s) dtds-\\
&\hspace{-5cm}-\iint_{\R^2}
\left\{\int_{\R}\Big\{e^{i(t-s)u}-1-
\frac{i(t-s)u}{u^2+1}\Big\}
\frac{d\sigma(u)}{u^2}\right\}
\frac{\partial^2}{\partial t\partial s}
\phi(t,s)dtds.
\end{split}
\]

The first integral on the right hand-side
vanishes since $\phi$ is a Schwartz function.
The function
\[
K(t,s,u)=\begin{cases}
\Big\{e^{i(t-s)}-1-\frac{i(t-s) u}{u^2+1}\Big\}
\frac{1}{u^2},~~{\rm if}~~u\neq 0,\\
~~~\\
-\frac{(t-s)^2}{2},~~~~~~~~~~~~~~~~~~~~~~~ {\rm
if}~~u=0,
\end{cases}
\]
is continuous since
\[
\lim_{u\to 0} \Big\{e^{i(t-s)}-1- \frac{i(t-s)
u}{u^2+1}\Big\}\frac{1}{u^2}=- \frac{(t-s)^2}{2}.
\]
Moreover, we have the bounds
\[
\left|K(t,s,u)\right|\leq\begin{cases}
\frac{(t-s)^2+\left|t-s\right|}
{u^2+1}~~{\rm if}~~\left|u\right|< 1,\\
~~~\\
\frac{4+\left|t-s\right|}{u^2+1}~~~~~~ {\rm
if}~~\left|u\right|\geq 1.
\end{cases}
\]
Therefore,
\[
\begin{split}
&\iint_{\R^2}\int_\R
\left|K(t,s,u)\frac{\partial^2}{\partial t\partial s}
\phi(t,s)\right|d\sigma(u)dtds\\
&\leq \iint_{\R^2}
\left|\frac{\partial^2}{\partial t\partial s}
\phi(t,s)\right|
\left\{\int_\R \left|K(t,s,u)
\right|d\sigma(u)\right\}dtds\\
&\leq
\iint_{\R^2}\left|\frac{\partial^2}{\partial
t\partial
s}\phi(t,s)\right|\left\{\int_{\left|u\right|<1}
\frac{(t-s)^2
+\left|t-s\right|}{u^2+1}d\sigma(u)\right\}dtds\\
&\hspace{10mm}+\iint_{\R^2}\left|
\frac{\partial^2}{\partial t\partial s}
\phi(t,s)\right|\left\{\int_{\left|u\right|\geq 1}
\frac{4+\left|t-s\right|}{u^2+1}
d\sigma(u)\right\}dtds\\
&<K\iint_{\R^2}\left|\frac{\partial^2}
{\partial t\partial s}\phi(t,s)
(\left|t-s\right|+2)^2\right|dtds<\infty\\
\end{split}
\]
where
\[
K=\int_\R \frac{d\sigma (u)}{u^2+1}<\infty.
\]
By Fubini's theorem and integration by parts we
obtain
\[
\begin{split}
\int_{\R^2}\left\{\int_{\R}\Big\{e^{i(t-s)u}-1-
\frac{i(t-s)u}{u^2+1}\Big\}\frac{d\sigma(u)}{u^2}
\right\}\frac{\partial^2}{\partial t\partial s}
\phi(t,s)dtds&=\\
&\hspace{-10cm}
=\int_{\R}\left\{\int_{\R^2}\left\{e^{i(t-s)u}-
1-\frac{i(t-s)u}{u^2+1}
\right\}\frac{1}{u^2}\frac{\partial^2}{\partial
t\partial s}\phi(t,s)dtds\right\}d\sigma(u)\\
&\hspace{-10cm} =
\int_{\R}\left\{\int_{\R^2}e^{i(t-s)u}
\phi(t,s)dtds\right\}d\sigma(u)\\
&\hspace{-10cm} =
\int_{\R^2}\left\{\int_{\R}e^{-i(s-t)u}d\sigma(u)\right\}\phi(t,s)dtds,
\end{split}
\]
and by \eqref{f_t} we conclude that
\[
\left\langle r''(t-s), \phi(t,s)\right\rangle=
\left\langle
\widehat{d\sigma}(s-t),\phi(t,s)\right\rangle.
\]
Thus,
\[
\frac{\partial^2}{\partial t\partial s}
K_r(t,s)=r''(t-s)=\widehat{d\sigma}(s-t).
\]
\mbox{}\qed\mbox{}\\

\section{The operator $T_m$}
\setcounter{equation}{0}
\label{tmmm}
We now focus on the case $d\sigma(u)=m(u)du$ in \eqref{kipour_2007}, where
$m$ is a positive and measurable function such that
\begin{equation}
\label{u2}
\int_{\R}\frac{m(u)du}{u^2+1}<\infty.
\end{equation}

We define an (unbounded in general) operator
$T_m$ by
\begin{equation}
\widehat{T_mf}(u) {\stackrel{{\rm
def.}}{=}}\sqrt{m(u)}\widehat{f}(u),
\label{Tm}
\end{equation}
where $\widehat{f}$ denotes the Fourier
transform of $f$; see \eqref{fourier}. The domain of
$T_m$,
\begin{equation}
\label{voltaire}
{\rm dom}~(T_m){\stackrel{{\rm def.}}{=}}
\left\{f\in
{\mathbf L}_2(\R):~~\int_{\R}m(u)
\left|\widehat{f}(u)
\right|^2du<\infty\right\},
\end{equation}
contains in particular the Schwartz space
$\s$ since $m$
satisfies \eqref{u2} and since the Fourier
transform maps $\s$ into itself.\\

When $m$ is summable, the integral in \eqref{voltaire} can be
rewritten as a double integral as explained in the previous
section:
\begin{equation}
\int_{\mathbb R}
m(u)|\widehat{f}(u)|^2du=
\iint_{\R^2}f(t){f(s)^*}\widehat{m}(t-s)dtds.
\label{phi}
\end{equation}

When
\begin{equation}
\label{mh}
m(u)=\frac{1}{2\pi}|u|^{1-2H},
\end{equation}
the operator $T_m$ reduces, up to a multiplicative constant, to
the operator $M_H$ defined in \cite[(2.10) p. 304]{eh} and in
\cite[ Definition 3.1 p. 354]{bosw}, and the function $r(t)$ in
\eqref{kipour_2007} is given by \eqref{mr}. We note that the set
\eqref{voltaire} has been introduced in \cite[Theorem 3.1, p.
258]{MR1790083} for $m$ of the form \eqref{mh}. Multiplying
\eqref{mh} by
\[
\frac{2\pi H(1-2H)}{\Gamma(2-2H)\cos(\pi H)},\] that is,
considering
\[
m(u)= \frac{H(1-2H)}{\Gamma(2-2H)\cos(\pi H)}
\left|u\right|^{1-2H}
\]
and using \cite[Formula 12, p. 170]{GS}, that is, in the sense of
distributions,
\[
\int_{\mathbb R}|x|^\lambda
e^{iux}dx=-2\Gamma(1+\lambda)\sin\left(\frac{\pi\lambda}{2}\right)
|u|^{-\lambda-1},\quad\lambda\not\in{\mathbb Z},
\]
leads to (with $\lambda=1-2H$)
\[
\begin{split}
\widehat{m}(u)&=-\frac{H(1-2H)\Gamma(2-2H)\sin\left(\frac{\pi(1-2H)}
{2}\right)}
{\Gamma(2-2H)\cos(\pi H)}\left|u\right|^{2H-2}\\
&= 2H(2H-1) \left|u\right|^{2H-2}
\end{split}
\]
See \cite[(2.1), p. 584]{MR1741154}. The norm
$\left|f\right|_{\phi}^2$ defined in \cite[(2.2),
p. 584]{MR1741154} is then equal to \eqref{phi}.\\

The operator $T_m$ will play a key role in the
sequel. We now study its main properties.\\

\begin{La} For every $t\in{\mathbb R}$,
the function
\[
I_t{\stackrel{{\rm def.}}{=}}
\begin{cases}1_{[0,t]},\quad{if}\quad t> 0,\\
                1_{[t,0]},\quad{if}\quad t< 0,
\end{cases}
\]
belongs to the domain of $T_m$.
\end{La}
{\bf Proof:} We consider the case $t> 0$. The
case where $t<0$ is treated in a similar way. We
have
\[
\begin{split}
\int_{\R}m(u)|\widehat{1_{[0,t]}}(u)|^2du&=
\int_{-\infty}^{-1}m(u)
\Big|\frac{e^{-itu}-1}{-iu}\Big|^2du+\\
&\hspace{5mm}+ \int_{-1}^1m(u)\Big|
\frac{e^{-itu}-1}{-iu}\Big|^2du+\\
&\hspace{5mm}+\int_1^{\infty}m(u)
\Big|\frac{e^{-itu}-1}{-iu}\Big|^2du.
\end{split}
\]
The first and last integrals converge in view of
\eqref{u2}, and the second is trivially
convergent.
\mbox{}\qed\mbox{}\\

\begin{Tm}
The operator $T_m$ is self-adjoint and closed.
It is bounded if and only if $m$ is bounded.
\end{Tm}
{\bf Proof:} For $f$ and $g$ in the domain of
$T_m$ we have
\[
\left\langle f,T_mg\right\rangle_{\L_2(\R)}
=\left\langle T_mf,g\right\rangle_{\L_2(\R)}.
\]
Thus, $T_m\subset T_m^*$ and the operator $T_m$
is hermitian. We show that it is self-adjoint:
let $g\in\Doms$. The map $f\to\left\langle
T_mf,g\right\rangle_{\L_2(\R)}$ is continuous and so is the
map
\[
f\to\left\langle
\widehat{T_mf},\widehat{g}\right\rangle_{\L_2(\R)}
=\left\langle \sqrt{m}\widehat{f},
\widehat{g}\right\rangle_{\L_2(\R)},
\]
and the map
\[
\widehat{f}\to\left\langle
\widehat{f},\sqrt{m}\widehat{g}\right\rangle_{\L_2(\R)}
\]
is also continuous. By Riesz representation
theorem, $\sqrt{m}\widehat{g}\in\L_2(\R)$, hence
$g\in\Dom$ and we get $T_m^*\subset T_m$.\\

We now show that $T_m$ is closed: let $f_n\rightarrow f$ and
$T_mf_n\rightarrow g$. We have $\widehat{f_n}\rightarrow
\widehat{f}$. Thus $T_mf_n\rightarrow g$ leads to
$\widehat{T_mf_n}\rightarrow\widehat{g}$, and thus
$\sqrt{m}\widehat{f_n}\rightarrow\widehat{g}$. By
\cite[Th\'eor\`eme 2.3, p. 95]{descombes} there exists a
subsequence $n_k$ such that
$\widehat{f_{n_k}}\rightarrow\widehat{f}$ point-wise a.e. and so
$\widehat{T_mf_{n_k}}\rightarrow\widehat{g}$ point-wise a.e., and
we have $\widehat{T_mf}=\widehat{g}$, a.e.
\\

Finally we show that the operator $T_m$ is bounded if and only if
$m$ is bounded. First, if $m$ is bounded then there exists a $K>0$
such that $\left|m(u)\right|<K$ for any $u\in\R$ and we get, for
any $f\in\L_2(\R)$
    \[
\int_{\R}m(u)|\widehat{f}(u)|^2du<
K\int_{\mathbb R} |\widehat{f}(u)|^2du,
    \]
    so $T_m$ is bounded since the Fourier transform is an isometry.
Now if $T_m$ is bounded, then there exists a $K\in\R$ such that
for any $f\in\L_2(\R)$
    \[
    \int_{\R}m(u)\left|\widehat{f}(u)
    \right|^2du\leq K\int_{\R}\left|\widehat{f}(u)
    \right|^2du.
    \]
Assume that $m$ is unbounded. Then for any $N\in\N$ there
     exists a measurable set $E_N$ such that
     $\lambda(E_N)>0$ ($\lambda$ denotes the Lebesgue measure)
and $m(u)\geq N$ on $E_N$, where, without loss of generality,
one may take $E_N$ such that
     $\lambda(E_N)\leq 1$.
Define $f_n$ such that
$\widehat{f_n}=1_{E_N}$ on $E_N$. We then have
    \[
  N\int_{E_N}\left|\widehat{f_n}(u)
  \right|^2du\leq
  \int_{E_N} m(u)
  \left|\widehat{f_n}(u)
  \right|^2du\leq K\int_{E_N}
  \left|\widehat{f_n}(u)\right|^2du,
    \]
hence $N\leq K$, but this is impossible,
so $m$ is bounded.
\mbox{}\qed\mbox{}\\

For $m(u)=\frac{1}{2\pi}|u|^{1-2H}$, we have that
\begin{equation}
\label{support}
{\rm supp}~T_mI_t\subset{\rm supp}~I_t,\quad
t\in{\mathbb R}.
\end{equation}
In general this property will not hold, as we now
illustrate with a
counterexample. This example is of particular
importance to mark the
difference between  our approach and the
approach presented in \cite{Alos2001}.

\begin{Ex}
{\rm  Let $m(u)=u^4e^{-2u^2}$. We have
\[
\begin{split}
(T_m1_{[0,t]})(s)&=\frac{1}{2\pi} \int_{\mathbb
R} e^{isu}u^2e^{-u^2}
\cdot\frac{e^{-itu}-1}{-iu}du\\
&=-\frac{1}{2\pi i}\int_{\mathbb R}
e^{isu}ue^{-u^2}(e^{-itu}-1)du\\
&=-\frac{1}{2\pi i}\int_{\mathbb R}
ue^{i(s-t)u}e^{-u^2}du\\
&\hspace{5mm}+\frac{1}{2\pi i}
\int_{\mathbb R}ue^{isu}e^{-u^2}du\\
&=\Phi(s)-\Phi(s-t),
\end{split}
\]
where
\[
\Phi(s)=\frac{1}{2\pi i}\int_{\mathbb R}
ue^{isu}e^{-u^2}du.
\]
We have:
\[
\Phi(s)=\frac{1}{2\pi i} \int_{\mathbb R}
ue^{isu}e^{-u^2}du=\frac{1}{2\pi
i}e^{-\frac{s^2}{4}}\int_{\mathbb R}
ue^{-(u-\frac{is}{2})^2}du
=\frac{s}{4\sqrt{\pi}}e^{\frac{-s^2}{4}}.
\]
Thus,
\[
(T_m1_{[0,t]})(s)= \frac{1}{4\sqrt{\pi}}
\left\{(t-s)e^{-\frac{(t-s)^2}{4}}+se^{-\frac{s^2}{4}}\right\}.
\]
The support of the function $T_m(I_t)$ is not bounded, and in
particular \eqref{support} is not in force.}
\label{Ex:u4}
\end{Ex}

When $m$ is bounded, we note that $T_m$ is a
translation invariant operator.\\

We now recall the definitions of the Hermite polynomials and of
the Hermite functions. Then, in Proposition \ref{pn:herm} below we
study the action of the operator $T_m$ on Hermite functions.
\begin{Dn} The Hermite polynomials
$\left\{h_n(x)\,\, n\in{\mathbb N}_0\right\}$ are defined by
\begin{equation*}
h_n(x){\stackrel{{\rm def.}}{=}}(-1)^n
e^\frac{x^2}{2}\frac{d^n}{dx^n}
(e^{-\frac{x^2}{2}}),~~~n=0,1,2\ldots.
\end{equation*}
\end{Dn}
\begin{Dn}
The Hermite functions are defined by
\begin{equation*}
\widetilde{h}_n(x){\stackrel{{\rm def.}}{=}}
\frac{h_{n-1}(\sqrt{2}x) e^{-\frac{x^2}{2}}}{\pi^{\frac{1}{4}}
\sqrt{(n-1)!}},~~~n=1,2,\ldots.
\end{equation*}
\end{Dn}
The following proposition the main properties of
the Hermite functions which we will need; see
\cite[p. 349]{bosw} and the references therein.
\begin{Pn} \cite[p. 349]{bosw}
The Hermite functions $\{\widetilde{h}_n,\,\,n\in{\mathbb
N}\}$ form
an orthonormal basis of ${\mathbf L}_2(\R)$.
Furthermore,
\begin{equation}
|\widetilde{h}_n(u)|\leq
\begin{cases}
Cn^{-\frac{1}{12}}~\quad{if}\quad|u|\leq 2\sqrt{n},\\
Ce^{-\gamma u^2}\quad{if}\quad|u|>2\sqrt{n},
\end{cases}
\label{mu}
\end{equation}
where $C$ and $\gamma>0$ are constants independent
of $n$. Finally, the Fourier transform of the
Hermite function is given by
\begin{equation}
\widehat{\widetilde{h}}_n(u)=
\sqrt{2\pi}(-1)^{n-1}\widetilde{h}_n(u).
\label{hermiteff}
\end{equation}
\label{pn:herm}
\end{Pn}
Using the previous proposition, we now study the
functions $T_m\widetilde{h}_n$.
\begin{Pn}
Assume that the function $m$ satisfies a bound of
the type:
\begin{equation}
m(u)\leq
\begin{cases}
K\left|u\right|^{-b}\quad{if}\quad|u|\leq 1,\\
K'~~~~~~\hspace{7mm}\quad{if}\quad|u|>1,
\end{cases}
\label{mbound}
\end{equation}
where $b<2$ and $0<K,K'<\infty$. Then,
\begin{equation}
\left|(T_m\widetilde{h}_n)(u)
\right|\leq\tilde{C_1}n^{\frac{5}{12}}
+\tilde{C_2},
\label{bound}
\end{equation}
where $\tilde{C_1}$ and $\tilde{C_2}$ are
constants independent of $n$. \label{T_mI_t}
\end{Pn}
{\bf Proof:} Using \eqref{hermiteff} we have
\[
\begin{split}
\left|(T_m\widetilde{h}_n)(u)\right|&=
\frac{1}{2\pi}\left|\int_\R e^{iuy}
\widehat{\widetilde{h}}_n(y)\sqrt{m(y)}dy\right|\\
&\leq\frac{1}{\sqrt{2\pi}} \int_\R
\left|\widetilde{h}_n(y) \right|\sqrt{m(y)}dy.
\end{split}
\]
We now compute an upper bound for the integral
$$\int_\R \left|\widetilde{h}_n(y)
\right|\sqrt{m(y)}dy=I_1+I_2+I_3,$$
where
\[
\begin{split}
&I_1=\int_{-\infty}^{-2\sqrt{n}}\left|
\widetilde{h}_n(y)\right|\sqrt{m(y)}dy,\\
&I_2=\int_{-2\sqrt{n}}^{2\sqrt{n}}\left|
\widetilde{h}_n(y)\right|\sqrt{m(y)}dy,\\
&I_3=\int_{2\sqrt{n}}^{\infty}\left|
\widetilde{h}_n(y)\right|\sqrt{m(y)}dy.
\end{split}
\]
By \eqref{mu} we have
\[
\begin{split}
I_1&\leq C\int_{-\infty}^{-2\sqrt{n}}e^{-\gamma
y^2} \sqrt{m(y)}dy,\\
I_2&\leq Cn^{-\frac{1}{12}}
\int_{-2\sqrt{n}}^{2\sqrt{n}}\sqrt{m(y)}dy,\\
I_3&\leq C\int_{2\sqrt{n}}^{\infty}e^{-\gamma
y^2}\sqrt{m(y)}dy.
\end{split}
\]
We have:
\[
\begin{split}
\int_{-2\sqrt{n}}^{2\sqrt{n}}\sqrt{m(y)}dy&=
\int_{-2\sqrt{n}}^{-1}
\sqrt{m(y)}dy+\int_{-1}^{1}\sqrt{m(y)}dy+
\\
&\hspace{5mm}+\int_{1}^{2\sqrt{n}}\sqrt{m(y)}dy\\
&\leq \sqrt{K'}\int_{-2\sqrt{n}}^{-1}dy+
\sqrt{K}\int_{-1}^{1}
\left|y\right|^{-\frac{b}{2}}dy+\\
&\hspace{5mm}+\sqrt{K'}\int_{1}^{2\sqrt{n}}dy\\
&=
2\sqrt{K'}(2\sqrt{n}-1)+4\sqrt{K}\frac{1}{2-b},
\end{split}
\]
so that we get
\[
I_2\leq C_1n^{\frac{5}{12}}+C_2n^{-\frac{1}{12}}.
\]
Furthermore,
\[
\begin{split}
\int_{2\sqrt{n}}^{\infty}e^{-\gamma
y^2}\sqrt{m(y)}dy &\leq
\sqrt{K'}\int_{2\sqrt{n}}^{\infty}e^{-\gamma
y^2}dy\\
&\le \sqrt{K'}\int_{\mathbb R}e^{-\gamma y^2}dy\\
&= 2\sqrt{K'}\sqrt{\frac{\pi}{\gamma}}=C_3.
\end{split}
\]
Finally we get
\[
\begin{split}
\left|(T_m\widetilde{h}_n)(y)\right|&
\leq\frac{1}{\sqrt{2\pi}}( C_1n^{\frac{5}{12}}
+C_2n^{-\frac{1}{12}}+2C_3)\\
&\le \tilde{C_1}n^{\frac{5}{12}}+
\tilde{C_2}
\end{split}
\]
for appropriate $\tilde{C_1}$ and $\tilde{C_2}$.
\mbox{}\qed\mbox{}\\

\begin{La}
The function $T_m \widetilde{h}_k$ is uniformly
continuous for every $k\in\N$. More precisely, it holds that
\begin{equation}
\label{eq:uni_conv}
\left|(T_m\widetilde{h}_k)(t)- (T_m\widetilde{h}_k)(s)\right|\leq
\left|t-s \right|\left\{ \widetilde{C}_1k^\frac{11}{12}
+\widetilde{C}_2\right\},
\end{equation}
\label{edf}
where $\widetilde{C}_1$ and $\widetilde{C}_2$ are constant independent
of $k$.
\end{La}
{\bf Proof:} Let $t,s\in\R$. We have
\[
(T_m\widetilde{h}_k)
(t)-(T_m\widetilde{h}_k)(s)
=\frac{1}{\sqrt{2\pi}}
(-1)^{k-1}\int_\R
 \left(e^{-iut}-e^{-ius}\right)\sqrt{m(u)}
\widetilde{h}_k(u)du.
\]
Taking into account that
\[
\left|e^{-iut}-e^{-ius}\right|\leq\left|u(t-s)\right|,
\]
we get
\[
\left|(T_m\widetilde{h}_k)
(t)-(T_m\widetilde{h}_k)(s)
\right|\leq \frac{\left|t-s\right|}{\sqrt{2\pi}}
\int_\R \left|u\right|\sqrt{m(u)}
\left|\widetilde{h}_k(u)\right|du.
\]
To conclude the proof, it suffices to show that $\int_\R
\left|u\right|\sqrt{m(u)}
\left|\widetilde{h}_k(u)\right|du<\infty$. By \eqref{mbound} and
\eqref{mu} we have
\begin{equation*}
\begin{split}
\int_\R \left|u\right|\sqrt{m(u)}
\left|\widetilde{h}_k(u)\right|du&\leq
Ak^{-\frac{1}{12}}
 \int_{\left|u\right|\leq 1}\left|u
 \right|^{1-\frac{b}{2}}du+\\
 &\hspace{5mm}+Bk^{-\frac{1}{12}}
 \int_{1<\left|u\right|\leq 2\sqrt{k}}\left|u
 \right|du+\\
&\hspace{5mm}
 +C\int_{\left|u\right|>2\sqrt{k}}
 \left|u\right|e^{-\gamma u^2}du\\
&\leq
\widetilde{C}_1k^\frac{11}{12}
+\widetilde{C}_2,
\end{split}
\end{equation*}
 where all the constants are independent of $k$.
\mbox{}\qed\mbox{}\\

We conclude this section with a remark.
\begin{Rk}
When
\[
\int_{\mathbb R}\frac{\ln m(u)}{u^2+1}du>-\infty,
\]
the function $m$ admits a factorization
$m(u)=|h(u)|^2$, where $h$ is an outer function.
One can define an operator $\widetilde{T}_m$ through $\widetilde{T}_mf
{\stackrel{{\rm def.}}{=}} \widehat{h}*f$ rather
than the operator $T_m$. We will not pursue this
direction here.
\end{Rk}

\section{The white noise space and the Brownian motion}
\label{S4}
\setcounter{equation}{0}
In this section we review the construction of the white noise
space and recall some results related to the Brownian motion. We
refer the reader to \cite{MR0451429}, \cite{MR1408433},
\cite{MR1387829} and \cite{MR1851117} for additional information
and references. To build the white noise space one considers the
subspace ${\mathcal S}_{\mathbb R}({\mathbb R})$ of the Schwartz
space which consists of {\sl real valued} functions. Denote
by ${\mathcal S}_{\mathbb R}({\mathbb R})^\prime$ its dual. Let
$\F$ be the $\sigma$-algebra of Borel sets in the space
${\mathcal S}_{\mathbb R}({\mathbb R})^\prime$. The function
\begin{equation*}
K(s_1-s_2)=
\exp(-\left\|s_1-s_2\right\|_{\L_2(\R)}^2/2)
\end{equation*}
is positive in  ${\mathcal S}_{\mathbb R}({\mathbb R})$ in the
sense of reproducing kernels since
\[
\begin{split}
\exp(-\left\|s_1-s_2\right\|_{\L_2(\R)}^2/2)&=
\exp(-\left\|s_1\right\|_{\L_2(\R)}^2/2)\times\\
&\hspace{5mm}\times \exp\langle s_1\, ,\,
s_2\rangle_{\L_2(\R)}\times
\exp(-\left\|s_2\right\|_{\L_2(\R)}^2/2).
\end{split}
\]
The space  ${\mathcal S}_{\mathbb R}({\mathbb R})$ is nuclear,
and therefore the Bochner--Minlos theorem (see for instance
\cite[Appendix A, p. 193]{MR1408433}) implies that there exists a
probability measure $P$ on $({\mathcal S}_{\mathbb R}({\mathbb
R})^\prime,\F)$ such that, for all $s\in{\mathcal S}_{\mathbb
R}({\mathbb R})$,
\begin{equation}
E(e^{iQ_s(s')}){\stackrel{{\rm def.}}{=}}
\int_{\ss}e^{iQ_s(s')}dP(s')
=e^{-\frac{\left\|s\right\|_{\L_2(\R)}^2}{2}},
\label{exp}
\end{equation}
where $Q_s$ denotes the linear functional $Q_s(s')=\left\langle
s',s\right \rangle_{{\mathcal S}_{\mathbb R}({\mathbb R})^\prime,
{\mathcal S}_{\mathbb R}({\mathbb R})}$; see \cite[(2.1) p.
348]{bosw}, \cite[(2.3) p. 303]{eh}. Note that $Q_s$ is the
canonical isomorphism of the Schwartz space ${\mathcal
S}_{\mathbb R}({\mathbb R})$ onto its bidual; see \cite[p.
7]{MR99f:60082}. Definition \eqref{exp} implies in particular that
\begin{equation}
E(Q_s)=0\quad{\rm and}\quad
E(Q_s^2)=\left\|s\right\|_{\L_2(\R)}^2.
\label{EQ}
\end{equation}
In view of \eqref{EQ}, the map
\begin{equation}
\label{eq:iso}
s\to Q_s
\end{equation}
is an isometry from the real Hilbert space ${\mathcal S}_{\mathbb
R}({\mathbb R}) \subset\L_2(\R)$ into the real Hilbert space
$\L_2({\mathcal S}_{\mathbb R}({\mathbb R})^\prime, \F,dP)$. It
extends to an isometry from $\L_2(\R)$ into $\L_2({\mathcal
S}_{\mathbb R}({\mathbb R})^\prime,\F,dP)$, and we define for
$f\in\L_2(\R)$
\begin{equation}
Q_f(s^\prime)\Def\lim_{n\to\infty}Q_{f_n}(s^\prime),
\end{equation}
where the limit is in $\L_2({\mathcal S}_{\mathbb R}({\mathbb
R})^\prime,\F,dP)$ and where $f_n\to f$ in $\L_2(\R)$. The limit
is easily
shown not to depend on $\left(f_n\right)$.\\

In the sequel we consider complex-valued functions. The map
\eqref{eq:iso} extends to an isometry between the complexified
spaces of $\L_2(\R)$ and $\L_2(\ss,\F,dP)$. See for instance
\cite[pp. V4-V5]{Bourbaki81EVT} for the
complexification of Hilbert spaces.\\

The triplet $\L_2({\mathcal S}_{\mathbb R}({\mathbb
R})^\prime,\F,P)$ is called the white noise space. In accordance
with the notation standard in probability theory, we set
\[
\Omega={\mathcal S}_{\mathbb R}({\mathbb R})^\prime,
\]
and denote by
\begin{equation*}
\W=\L_2(\Omega,\F,P).
\end{equation*}
the complexified space of $\L_2({\mathcal S}_{\mathbb R}({\mathbb
R})^\prime,\F,P)$.\\

The Brownian motion is a family $\{B(t,\omega)\}$ of random
variables in the white noise space with the following property:
\begin{enumerate}
    \item $B(0,\omega)=0$ almost surely with
    respect to $P$.
    \item $\{B(t,\omega)\}$ is a Gaussian stochastic
    process with mean zero and $B(t,\omega)$ and
    $B(s,\omega)$ have the covariance $\rm{min}(t,s)$.
    \item $s\to B(s,\omega)$ is a continuous for
    almost all $\omega$ with respect to $P$.
\end{enumerate}
Define the stochastic process
\begin{equation*}
\widetilde{B}(t,\omega)=Q_{I_t}(\omega),\quad t\in{\mathbb R}.
\end{equation*}
Then, for $t,s\ge 0$,

\[
E(\widetilde{B}(t,\omega)\widetilde{B}(s,\omega)^*)=\left\langle I_t,
I_s\right\rangle_{\L_2(\R)}=\int_\R
I_t(u){(I_s(u))^*}du=\min(t,s).
\]

By Kolmogorov's continuity Theorem the process
$\left\{\widetilde{B}(t,\omega)\right\}$ has a continuous version
$\left\{B(t,\omega)\right\}$, with is a Brownian motion. For
$F\in\W$ we now recall Wiener-Ito chaos expansion. In the
stochastic process literature this expansion is for real valued
functions. We write it for the complexification of the underlying
Hilbert spaces. This creates no technical problem.\\

The white noise probability space $\W$ admits a special
orthonormal basis $\{H_{\alpha}\}$, indexed by the set $\ell$
(defined by \eqref{eq:ell}). We will not recall the definition of
this basis here, and refer to \cite[Definition 2.2.1 p.
19]{MR1408433}.

\begin{Pn}(Wiener-Ito chaos expansion)
Every $F\in\W$ can be
written as
\begin{equation*}
F=\sum_{\alpha\in\ell}c_{\alpha}H_{\alpha},
\end{equation*}
with $\alpha\in\ell$, $c_\alpha\in\C$,
 and
\[
\|F\|_{\W}^2=
\sum_{\alpha\in\ell}\alpha!\left| c_{\alpha}\right|^2<\infty,
\]
 where $\alpha!=\alpha_1 !\alpha_2! \alpha_3 !\cdots$ and
\begin{equation}
H_\alpha(\omega){\stackrel{{\rm def.}}{=}}
\prod_{k=1}^\infty h_{\alpha_k}
\Big(Q_{\widetilde{h}_k}(\omega)\Big),\quad\omega\in\Omega.
\label{H}
\end{equation}
\label{pn:wh}
\end{Pn}

\section{The Kondratiev space and the Wick
product} \setcounter{equation}{0}
\label{S5}
The Wick product is defined through:
\begin{Dn}
Let $\alpha,\beta\in\ell$, then
\begin{equation*}
H_\alpha\diamond H_\beta= H_{\alpha+\beta}.
\end{equation*}
\end{Dn}
\begin{Dn}
\label{Wick}
Let $F,G$ be two elements in $\W$
\[
F=\sum_{\alpha}a_{\alpha\in\ell}H_{\alpha},\quad and\quad
G=\sum_{\alpha\in\ell}b_{\alpha}H_{\alpha},
\]
where $\alpha\in\ell$, $a_\alpha,b_\alpha\in\C$ and $a_\alpha,
b_\alpha\neq 0$ for only a finite number of indexes $\alpha$. The
Wick product of $F$ and $G$ is defined by
\[
(F\diamond G)(\omega)=\sum_{\alpha,\beta\in\ell}
a_{\alpha}b_{\beta}H_{\alpha+\beta}(\omega)=
\sum_{\gamma}\left(\sum_{\gamma=\alpha+\beta}
a_{\alpha}b_{\beta}\right)H_{\gamma}(\omega).
\]
\end{Dn}
This product can be shown to be independent of the basis
$\left\{{H_\alpha}\right\}_{\alpha\in\ell}$, see
\cite[Appendix D, p. 209]{MR1408433},
but it is not defined for all pairs of elements in the white
noise space. See \cite{MR1408433}.\\

The Kondratiev space $S_{-1}$ seems to be the most convenient
space within which the Wick product is well defined. It is a
space of distributions. We first recall on which space of test
functions its elements operate.
\begin{Dn}
The Kondratiev space $S_1$ of stochastic test
functions consists of the elements in the form $f=\sum_{\alpha\in\ell}
a_\alpha H_\alpha\in{\mathcal W}$ such that
\begin{equation*}
\sum_{\alpha\in\ell} \left|a_\alpha\right|^2(\alpha !)^2
(2\N)^{k\alpha}<\infty, \quad k=1,2,\ldots,
\end{equation*}
and where
\[
(2{\mathbb
  N})^\alpha\stackrel{\rm def.}{=}
2^{\alpha_1}(2\cdot 2)^{\alpha_2}(2\cdot3)^{\alpha_3}
  \cdots,\quad
\alpha\in\ell.
\]
\end{Dn}

\begin{Dn}
The Kondratiev space $S_{-1}$ of stochastic
distributions consists of the elements in the form $F=\sum_{\alpha\in\ell}
b_\alpha H_\alpha$ with the property that
\begin{equation*}
\sum_{\alpha\in\ell}
\left|b_\alpha\right|^2(2\N)^{-q\alpha}<\infty,
\end{equation*}
for some $q\in\N$.
\end{Dn}

$S_{-1}$ can be identified with the dual of $S_1$
and the action of $F\in S_{-1}$ on $f=\sum_{\alpha\in\ell}a_\alpha
H_\alpha\in S_1$ is
given by
\begin{equation}
\left\langle F,f\right \rangle_{S_{-1},S_1}
\Def\sum_{\alpha\in\ell}\alpha!a_\alpha b_\alpha.
\label{action}
\end{equation}

We also note the following: let $\alpha\in\ell$. By \eqref{H},
using a Wick product calculation, we have
\begin{equation*}
H_\alpha (\omega)=\prod_{k=1}^\infty \left(Q_{\widetilde{h}_k}
(\omega)\right)^{\diamond\alpha_k} \label{}
\end{equation*}
for $\alpha=\epsilon^{(k)}=(0,0,
\cdots,0,1,0,\cdots),~\alpha_i=0$ for $i\neq k$ and $\alpha_k=1$
we get
\begin{equation*}
H_{\epsilon^{(k)}}=Q_{\widetilde{h}_k} =\int_\R
\widetilde{h}_k(t)dB(t).
\end{equation*}

We now review the main results associated with the
Wick product and the Hermite transform.\\

A key property of the basis $\{H_\alpha,\,\alpha\in\ell\}$ is the
following: define a map $\textbf{I}$ such that
$$ \textbf{I}(H_\alpha)=z^\alpha,$$
where $\alpha\in\ell$,
$z=(z_1,z_2,\ldots)\in\C ^{\N}$ (the set of
all sequences of complex numbers) and
$$z^\alpha=z_1^{\alpha_1}z_1^{\alpha_2}\cdots.$$
Then
$$\textbf{I}(H_\alpha \diamond H_\beta)=\textbf{I}
(H_\alpha)\textbf{I}(H_\beta).$$ The map
$\textbf{I}$ is called the Hermite transform.\\

We note that the spaces $S_{-1}$ and $S_1$ are
closed under the Wick product; see \cite[lemma 2.4.4, p 42]{MR1408433}.
\begin{Dn} Let $F=\sum_{\alpha\in\ell} a_\alpha
H_\alpha \in S_{-1}$. Then the Hermite transform of $F$, denoted
by $\textbf{I}(F)$ or $\widetilde{F}$, is defined by
$$\textbf{I} (F)(z)=\widetilde{F}(z)=
\sum_{\alpha\in\ell} a_\alpha z^\alpha.$$
\end{Dn}
\begin{Pn}\cite[Proposition 2.6.6, p. 59]{MR1408433}
\label{hermite product}
Let $F,G\in S_{-1}$. Then
$$\textbf{I}(F\diamond G)(z)=
(\textbf{I}(F)(z))\cdot(\textbf{I}(G)(z)).$$
\end{Pn}

\section{The $m$-Brownian motion
associated with $T_m$ } \setcounter{equation}{0}
In this section we define the $m$-Brownian
motion associated with the operator $T_m$ defined
in \eqref{Tm}, and the analogue of the
Wiener-Ito chaos expansion (see Proposition \ref{pn:wh}).
Using Lemma \ref{T_mI_t} we have
$T_mI_t\in\L_2(\R)$ and by the expansion in
$\L_2(\R)$ in term of the Hermite functions
$\widetilde{h}_n$ we obtain
\[
T_mI_t=\sum_{k=1}^\infty\left
\langle T_mI_t,
\widetilde{h}_k\right\rangle\widetilde{h}_k,
\]
where
\[
\left\langle T_mI_t,
\widetilde{h}_k\right\rangle=\left\langle
I_t,T_m\widetilde{h}_k\right\rangle=\int_\R
I_t(y)(T_m\widetilde{h}_k)(y)dy.
\]

\begin{Dn}
The process $\left\{ \widetilde{B}_m(t,\omega),\, t\in{\mathbb
R}\right\}$ defined through
\begin{equation*}
\widetilde{B}_m(t,\omega) {\stackrel{{\rm
def.}}{=}}Q_{T_mI_t}(\omega),
\end{equation*}
where $t\in\R$ and $\omega\in\Omega$ will be called the
$m$-Brownian motion associated with $m$.
\end{Dn}

As already noted in Section \ref{tmmm}, when
the function $m$ is given by \eqref{mh} (or, equivalently, the
function $r$ is given by \eqref{mr}), the operator $T_m$ is equal to
the operator $M_H$ defined in \cite{eh} and \cite{bosw}. Then, $B_m$
reduces to the fractional Brownian motion with
Hurst parameter $H\in (0,1)$.
\begin{La}
\label{expectation}
 The $m$-Brownian motion
has the following properties:
\begin{enumerate}
\item $E(\widetilde{B}_m(t,\omega)
\widetilde{B}_m(s,\omega)^*)=
K_r(s,t)$.
\item
\label{r1}$E\left(\left|\widetilde{B}_m(t,\omega)
-\widetilde{B}_m(s,\omega)
\right|^2\right)=2~{\rm Re}~r(t-s)$.
\item \label{r2}${\rm Re}~r(t)\leq C_1t^2+C_2t$
\end{enumerate}
for some positive constants $C_1$ and $C_2$.
\end{La}
{\bf Proof:} To prove item $(1)$ we note that
\[
\begin{split}
E(\widetilde{B}_m(t,\omega)
\widetilde{B}_m(s,\omega)^*)&=\left\langle
T_mI_t,T_mI_s
\right\rangle_{\L_2(\R)}\\
&=\int_{\R} (T_mI_t)(u)((T_mI_t)(u))^* du\\
&=
\int_{\mathbb R} \widehat{(T_mI_t)}(u)
({\widehat{(T_mI_s)}}(u))^*du\\
&=
\int_{\mathbb R}m(u)
\widehat{(1_{[0,t]})}(u)\left({
\widehat{(1_{[0,s]})}}(u)\right)^*du\\
&=
\int_{\mathbb R}m(u)\Big\{\int_{0}^t
e^{-iux}dx\Big\}\Big\{\int_{0}^se^{iuy}dy\Big\}du\\
&=
\int_{\mathbb R}
\frac{e^{-itu}-1}{u}\frac{e^{isu}-1}{u}m(u)du\\
&=
K_r(s,t).\\
\end{split}
\]
The proof of the second statement is carried out by direct
computations:
\[
\begin{split}
E&\left(\left|
\widetilde{B}_m(t,\omega)-
\widetilde{B}_m(s,\omega)\right|^2\right)\\
&=E((\widetilde{B}_m(t,\omega)-
\widetilde{B}_m(s,\omega))
(\widetilde{B}_m(t,\omega)
-\widetilde{B}_m(s,\omega))^*)\\
&=E\left\{\widetilde{B}_m(t,\omega)
\widetilde{B}_m(t,\omega)^*-
\widetilde{B}_m(t,\omega)
\widetilde{B}_m(s,\omega)^*\right.\\
&\left.\hspace{15mm}-\widetilde{B}_m(s,\omega)
\widetilde{B}_m(t,\omega)^*
+\widetilde{B}_m(s,\omega)
\widetilde{B}_m(s,\omega)^*\right\}.
\end{split}
\]
Using the first statement we get
\[
\begin{split}
&=E(\widetilde{B}_m(t,\omega)
\widetilde{B}_m(t,\omega)^*)-
E(\widetilde{B}_m(t,\omega)
\widetilde{B}_m(s,\omega)^*)\\
&\hspace{5mm}-E(\widetilde{B}_m(s,\omega)
\widetilde{B}_m(t,\omega)^*)+
E(\widetilde{B}_m(s,\omega)
\widetilde{B}_m(s,\omega)^*)\\
&=\left\{K_r(t,t)-K_r(s,t)-
K_r(t,s)+K_r(s,s)\right\}\\
&=(r(t)+r(t)^*-(r(s)+r(t)^*-r(s-t))\\
&\hspace{5mm}-(r(t)+r(s)^*-r(t-s))+r(s)+r(s)^*)\\
&=(r(t-s)+r(t-s)^*)\\
&=2~{\rm Re}~r(t-s).\\
\end{split}
\]
Finally, recall that we have
\[
{\rm Re}~r(t)=\int_{\R}
\Big\{1-\cos(tu)\Big\}\frac{m(u)}{u^2}du.
\]
Then, when $m(u)$ satisfies \eqref{bound} we get
\[
\begin{split}
{\rm Re}~r(t)&\leq 2\left\{K\int_0^1
\frac{\left|1-\cos(tu)\right|}{u^{2+b}}du
+K'\int_1^\infty
\frac{\left|1-\cos(tu)\right|}{u^2}du\right\}\\
&=
 2\left\{K
 \int_0^1\frac{2\sin^2(\frac{tu}{2})}{u^{2+b}}du
+K'\int_1^\infty
\frac{\left|1-\cos(tu)\right|}{u^2}du\right\}.
\end{split}
\]
But
\[
\int_0^1\frac{2\sin^2(\frac{tu}{2})}{u^{2+b}}
du\leq t^2\int_0^1\frac{u^2}{2u^{2+b}}du=
\frac{t^2}{2}\int_0^1\frac{1}{u^b}du=
\frac{t^2}{2(1-b)}
\]
since for $t\in [0,1]$ we get $tu\in [0,1]$ and
$\sin^2(\frac{tu}{2})\leq \frac{(tu)^2}{4}$.
Furthermore,
\[
\begin{split}
\int_1^\infty\frac{\left|1-
\cos(tu)\right|}{u^2}du &\leq
\int_0^\infty\frac{\left|1-
\cos(tu)\right|}{u^2}du\\
&=t\int_0^\infty
\frac{\left|1-\cos v\right|}{v^2}dv\\
&=t\int_0^1\frac{\left|1-\cos v\right|}{v^2}
dv+t\int_1^\infty\frac{\left|1-\cos v\right|}{v^2}
dv\\
&\leq t\int_0^1\frac{ 2\sin^2(\frac{v}{2})}{v^2}
dv+2t\int_1^\infty\frac{1}{v^2}dv\\
&\leq \frac{t}{2}\int_0^1\frac{v^2}{v^2}dv+2t
\int_1^\infty\frac{1}{v^2}dv\\
&=\frac{t}{2}+2t=\frac{5t}{2}.
\end{split}
\]
Thus,
\[
{\rm Re}~r(t)\leq\left|{\rm Re}~r(t)
\right|\leq 2\left\{K\frac{t^2}{2(1-b)}+K'
\frac{5t}{2}\right\}=C_1 t^2+C_2 t.
\]
\mbox{}\qed\mbox{}\\

The next step is to show that
$\left\{\widetilde{B}_m(t,\omega),\, t\in{\mathbb R}\right\}$
meets the criterion of Kolmogorov Theorem concerning the
existence of a continuous version of a given stochastic process.
Using the fact (see for instance \cite[p.5]{MR99f:60082} with
$p=2n$) that
\begin{equation}
E\Big(\left|\widetilde{B}_m(t,\omega)
\right|^{2n}\Big)=\kappa(2n)^{2n}
\Big(E\left(\left|\widetilde{B}_m(t,\omega)
\right|^2\right)\Big)^{\frac{2n}{2}}
\label{gaus}
\end{equation}
where
\[
\kappa(2n)=\sqrt{2}\left(
\frac{\Gamma(\frac{2n+1}{2})}{\sqrt{\pi}}
\right)^\frac{1}{2n}=\sqrt{2}
\left(\frac{2n!}{4^{n}n!}\right)^{\frac{1}{2n}}
\]
we have
\[
E\left(\left|
\widetilde{B}_m(t,\omega)
-\widetilde{B}_m(s,\omega)\right|^4\right)=
\kappa(4)^{4}E\Big(\left|
\widetilde{B}_m(t,\omega)
-\widetilde{B}_m(s,\omega)\right|^2\Big)^{2}.
\]
By \eqref{r1}, \eqref{r2} we get
\[
\kappa(4)^{4}E
\Big(\left|
\widetilde{B}_m(t,\omega)-\widetilde{B}_m(s,\omega)
\right|^2\Big)^{2}
=\kappa(4)^{4}\left(
\textsl{Re} ~r(t-s)\right)^2
\]
\[
\leq\kappa(4)^{4}\left(C_1 (t-s)^2+C_2 (t-s)
\right)^2=(t-s)^2(A+B(t-s))^2,
\]
for $t-s\in [0,1]$. By Kolmogorov's continuity theorem the
process $\widetilde{B}_m(t,\omega)$ has a continuous version
$B_m(t,\omega)$ where $t\in [0,1]$, which we will define as the
{\sl $m$-Brownian motion associated with $T_m$}. One can show in a
similar way that a continuous version exists on every finite
interval.

\begin{Pn}
$B_m(t)$ is a Gaussian random variable with
\[
E(B_m^n(t,\omega))=\begin{cases}
0,\quad\hspace{2.3cm}{\rm if}\quad n=2k-1\\
\dfrac{(2k)!}{2^k k!}\left\| T_mI_t\right\|^{2k},\quad{\rm
if}\quad n=2k
\end{cases}
\]
for $k=1,2,\ldots$.
\end{Pn}
{\bf Proof:} By \eqref{exp} with
$B_m(t,\omega)=Q_{T_mI_t}(\omega)$, we have
\[
E(\exp(i B_m(t,\omega)))=
e^{-\frac{ \left\|T_mI_t\right\|^2}{2}},
\]
and therefore $B_m(t,\omega)$ is a Gaussian
random variable. We now verify that
\[
E\Big(\sum_{n=0}^\infty \frac{i^n}{n!}B_m^n(t,\omega)\Big)
=\sum_{n=0}^\infty \frac{i^n}{n!}E(B_m^n(t,\omega)).
\]
Since $e^{-\frac{ \left\|T_mI_t\right\|^2}{2}}\in\R$, we have that:
$$E\Big(\sum_{n=0}^\infty \frac{i^n}{n!}B_m^n(t,\omega)\Big)=
E\Big(\sum_{n=0}^\infty
\frac{(-1)^{n}}{2n!}B_m^{2n}(t,\omega)\Big).$$
Since
$B_m(t,\omega)$ is a centered Gaussian random variable, we have
that $E(B_m^{2k-1}(t,\omega))=0$ for $k=1,2,\ldots$, and we get
$$\sum_{n=0}^\infty \frac{i^n}{n!}E(B_m^n(t,\omega))=
\sum_{n=0}^\infty \frac{(-1)^n}{2n!}E(B_m^{2n}(t,\omega)).$$
We have to verify
$$E\Big(\sum_{n=0}^\infty \frac{(-1)^n}{2n!}
B_m^{2n}(t,\omega)\Big)=
\sum_{n=0}^\infty \frac{(-1)^n}{2n!}E(B_m^{2n}(t,\omega)).$$
Let $\epsilon\in\R$.
Using \eqref{gaus} we have
\[
\begin{split}
\sum_{n=0}^\infty
\frac{\left|\epsilon
\right|^{2n} E(\left| B_m(t,\omega)
\right|^{2n})}{2n!}
&=\sum_{n=0}^\infty \frac{\left|\epsilon
\right|^{2n}
\kappa(2n)^{2n}
\left(E\left|B_m(t,\omega)\right|^2\right)^n}
{2n!}\\
&=\sum_{n=0}^\infty \frac{\left|
\epsilon\right|^{2n}
\left(E\left|B_m(t,\omega)\right|^2\right)^n}{2^nn!}
<\infty.
\end{split}
\]
We can thus use the dominated convergence
theorem, to obtain that
\[
\begin{split}
E\Big(\sum_{k=0}^\infty
\epsilon^k\frac{i^k}{k!}
B_m^k(t,\omega)\Big)&=
\sum_{k=0}^\infty \epsilon^k\frac{i^k}{k!}
E(B_m^k(t,\omega))\\
&=
\sum_{\ell =0}^\infty (-1)^\ell\epsilon^{2\ell}
\frac{\left\|T_mI_t\right\|^{2\ell}}{2^\ell \ell !}.
\end{split}
\]
The proof is completed by comparing the powers of
$\epsilon$ on both sides.
\mbox{}\qed\mbox{}\\

\begin{Rk}
In view of \eqref{qaz} we have
\[
\|T_mI_t\|^2=K_r(t,t).
\]
\end{Rk}

\begin{Rk}
Since for any $t\in{\mathbb R}$, $B_m(t)$ is written as a weighted sum
of the $\left\{H_\alpha,\,\alpha\in\ell\right\}$ (for an explicit
expression, see \eqref{BM} below), in turn being jointly Gaussian
random variables, it follows that $\left\{B_m(t),\,\ t\in{\mathbb
R}\right\}$ is a Gaussian process.
\end{Rk}
The following proposition will be used in the sequel of this
paper, where, as already noted, we develop the stochastic analysis
associated with the processes $B_m$.

\begin{Pn}
Let $f\in {\rm dom}~(T_m)$ and $n\in\N$. It holds that:
\[
Q_{T_mf}^{\diamond n}
(\omega)=n!\sum_{k=\left\lceil
\frac{n}{2}\right\rceil}^n
\Big(-\frac{1}{2}\Big)^{n-k}
\frac{Q_{T_mf}^{2k-n}(\omega)}{(2k-n)!}\frac{(\left\|
T_mf\right\|^2)^{n-k}}{(n-k)!}.
\]
In particular, for $f= I_t$, it holds that
\[
B_m^{\diamond n}(t)=n!\sum_{k=\left\lceil
\frac{n}{2}\right\rceil}^n \Big(-\frac{1}{2}
\Big)^{n-k}\frac{B_m^{2k-n}(t)}{(2k-n)!}
\frac{(\left\|T_mI_t\right\|^2)^{n-k}}{(n-k)!}.
\]
\end{Pn}
{\bf Proof:} Let $\epsilon\in\R$ then
\[
\exp^\diamond(Q_{\epsilon T_mf}(\omega))=
\sum_{n=0}^\infty \frac{(Q_{\epsilon T_mf}
(\omega))^{\diamond n}}{n!}=\sum_{n=0}^\infty
\frac{\epsilon^n(Q_{T_mf}(\omega))^{\diamond n}}{n!}.
\]
By \cite[Theorem 3.33, p. 32]{MR99f:60082}, we have
\[
\begin{split}
\exp^\diamond(Q_{\epsilon T_mf}
(\omega))&=\exp(Q_{\epsilon T_mf}(\omega)-
\frac{1}{2}\left\|\epsilon T_mf\right\|^2)\\
&=\sum_{k=0}^\infty \frac{(\epsilon
Q_{T_mf}(\omega)- \frac{1}{2}\epsilon^2
\left\|T_mf\right\|^2)^k}{k!}\\
&= \sum_{k=0}^\infty\sum_{j=0}^k \epsilon^{2k-j}
\frac{(Q_{T_mf} (\omega))^j}{j!}
\frac{(-\frac{1}{2}
\left\|T_mf\right\|^2)^{k-j}}{(k-j)!}.
\end{split}
\]
Hence,
\[
\sum_{n=0}^\infty
\frac{\epsilon^n(Q_{T_mf}(\omega))^{\diamond
n}}{n!} =\sum_{k=0}^\infty\sum_{j=0}^k
\epsilon^{2k-j}
\frac{(Q_{T_mf}(\omega))^j}{j!}\frac{(-\frac{1}{2}
\left\|T_mf\right\|^2)^{k-j}}{(k-j)!},
\]
and comparing the powers of $\epsilon$ leads to:
\[
\frac{(Q_{T_mf}(\omega))^{\diamond n}}{n!}
=\sum_{k=\left\lceil \frac{n}{2}\right\rceil}^n
\Big(-\frac{1}{2}\Big)^{n-k}\frac{(Q_{T_mf}
(\omega))^{2k-n}}{(2k-n)!}\frac{(\left\|
T_mf\right\|^2)^{n-k}}{(n-k)!}.
\]
\mbox{}\qed\mbox{}\\

\section{The $m$-white noise}
\setcounter{equation}{0}
One important aspect of the white noise space
theory is that the Brownian motion admits a
derivative, which belongs to the Hida space $(S)^*$ (the definition
of which we do not recall in this paper), and in particular in the
Kondratiev space $S_{-1}$. See
\cite[p. 53]{MR1408433}. In this section we prove that
this result still holds for the $m$-Brownian motion.
For the next definition, see also
\cite[Definition 2.5.5, p. 49]{MR1408433}, where the integral is
defined to be an element in the Hida space
$(S)^*$.
\begin{Dn}
Suppose that $Z:\R\to S_{-1}$ is a given function
with the property that
\begin{equation*}
\left\langle Z(t),f\right\rangle\in\L_1(\R,dt)
\end{equation*}
for all $f\in S_1$. Then $\int_\R Z(t)dt$ is
defined to be the unique element of $S_{-1}$ such
that
\begin{equation*}
\left\langle\int_\R Z(t)dt,f\right\rangle= \int_\R\left\langle
Z(t),f\right\rangle dt
\end{equation*}
for all $f\in S_1$.
\end{Dn}

In view of Lemma \ref{edf} the coefficients of the expansion
\eqref{white} below are continuous functions,
and not merely elements of ${\mathbf L}_2({\mathbb R})$.
\begin{Dn}
The $m$-white noise $W_m(t)$ is defined by
\begin{equation}
W_m(t)=\sum_{k=1}^\infty
\left(T_m\widetilde{h}_k\right)
(t)H_{\epsilon^{(k)}}.
\label{white}
\end{equation}
\end{Dn}
\begin{Tm} For every real $t$ we have that
$W_m(t)\in S_{-1}$, and it holds that
\begin{equation}
\label{bat_yam_290509}
B_m(t)=\int_0^t W_m(s)ds,
\quad  t\in{\mathbb R}.
\end{equation}
\end{Tm}
{\bf Proof:} Let $q\ge 2\in{\mathbb N}$. Then, using \eqref{bound},
we have:
\[
\sum_{k=1}^\infty \left| \left(T_m\widetilde{h}_k\right)(t)
\right|^2(2k)^{-q}\leq\sum_{k=1}^\infty
(\tilde{C_1}k^{\frac{5}{12}}+\tilde{C_2})^2(2k)^{-q} <\infty,
\]
and so $W_m(t)\in S_{-1}$. We now prove
\eqref{bat_yam_290509}. By construction, $B_m(t)\in{\mathcal W}$
for every $t\in{\mathbb R}$, and we can write
\begin{equation}
B_m(t)=\sum_{k=1}^\infty b_k(t)H_{\epsilon^{(k)}},
\label{BM}
\end{equation}
where
\[
 b_k(t)=\int_0^t
\left(T_m\widetilde{h}_k\right)(s)ds,
\]
with the convergence in the topology of ${\mathcal W}$.
We want to show that, for every $f\in S_1$, we have
\begin{equation*}
\langle B_m(t),f\rangle_{S_{-1},S_1}=\int_0^t\langle
W_m(u),f\rangle_{S_{-1},S_1}du,
\end{equation*}
where $\langle\,\cdot\, .\, \cdot\,\rangle_{S_{-1},S_1}$ denotes
the duality between $S_1$ and $S_{-1}$ (see \eqref{action}). To
that purpose, let $q\ge 2\in{\mathbb N}$. By using the estimate
\eqref{bound}, then, with $f=\sum_{\alpha\in\ell} f_\alpha
H_\alpha$ we have for $u\in[0,t]$ (and in fact for every $u\ge
0$),
\[
\begin{split}
\sum_{k=1}^\infty |T_m\widetilde{h}_k(u)f_k|&=\sum_{k=1}^\infty
|T_m\widetilde{h}_k(u)|(2k)^{-q}(2k)^q|f_k|\\
&\le\left(\sum_{k=1}^\infty(\widetilde{C}_1n^{\frac{5}{12}}+
\widetilde{C}_2)^2(2k)^{-2q}\right)^{\frac{1}{2}}\cdot
\left(\sum_{k=1}^\infty(2k)^{2q}|f_k|^2\right)^{\frac{1}{2}}\\
&<\infty,
\end{split}
\]
since $f\in S_1$. Therefore the series
\[
\sum_{k=1}^\infty |T_m\widetilde{h}_k(u)f_k|
\]
converges absolutely. Using the dominated convergence theorem we
can write
\[
\begin{split}
\int_0^t\langle W_m(u),f\rangle_{S_{-1},S_1}du&= \int_0^t
\left(\sum_{k=1}^\infty T_m\widetilde{h}_k(u)f_k\right)du\\
&=\sum_{k=1}^\infty \left(\int_0^tT_m\widetilde{h}_k(u)du\right)f_k\\
&=\langle B_m(t),f\rangle_{S_{-1},S_1}.
\end{split}
\]

\mbox{}\qed\mbox{}\\

We now show that, conversely,
\[
B_m(t)^\prime=W_m(t),
\]
in the sense of $S_{-1}$-processes; see \cite[p. 77]{MR1408433}
and further below in the current section. In the following
statements, the set $K_q(\delta)$ is defined by
$$ K_q(\delta)=
\{z\in\C^\N: \sum_{\alpha\in\ell} \left|z^\alpha\right|^2
(2\N)^{q\alpha}<\delta^2\}.$$ See \cite[Definition 2.6.4 p.
59]{MR1408433}.

\begin{Pn}
\label{pn:bd}
The function $\textbf{I}(W_m(t))(z)$ is bounded for
$(t,z)\in {\mathbb R}\times K_2(\delta)$.
\end{Pn}
{\bf Proof:} Write
\[
W_m(t)=\sum_{k=1}^\infty (T_m \widetilde{h}_k)(t)
Q_{\widetilde{h}_k}.
\]
Taking the Hermite transform we have
\[
\textbf{I}(W_m(t))(z)= \sum_{k=1}^\infty
(T_m\widetilde{h}_k)(t)z_k.
\]
Thus, for every $q\ge 2\in{\mathbb N}$  and using \eqref{bound} we
have: \vspace{-0cm}
\[
\begin{split}
\left|\textbf{I}(W_m(t))\right|
&=\left|\sum_{k=1}^\infty
(T_m\widetilde{h}_k)(t)z_k\right|\\
&=\left|\sum_{k=1}^\infty
(T_m\widetilde{h}_k)(t)
(2k)^\frac{q}{2}(2k)^
\frac{-q}{2}z_k\right|\\
&=\left(\sum_{k=1}^\infty \left|(T_m
\widetilde{h}_k)(t)\right|^2(2k)^q\right)^\frac{1}{2}
\left(\sum_{k=1}^\infty (2k)^{-q}\left|z^{\epsilon_k}
\right|^2\right)^\frac{1}{2}\\
&\leq \left(\sum_{k=1}^\infty \left\{\widetilde{C}_1k^\frac{5}{12}
+\widetilde{C}_2\right\}^2(2k)^q\right)^\frac{1}{2}
\left(\sum_{\alpha\in\ell} (2\N)^{-q\alpha}\left|z^{\alpha}
\right|^2\right)^\frac{1}{2}.
\end{split}
\]
The first sum converges when $q\ge 2$ and the second converges
since $z\in K_2(\delta)$. We conclude that the function
$\textbf{I}(W_m(t))(z)$ is bounded for any pair $(t,z)\in \R\times
K_2(\delta)$. \mbox{}\qed

\begin{Tm}
The function $\textbf{I}(W_m(t))(z)$ is
uniformly continuous in $t$ for $z\in K_4(\delta)$.
\end{Tm}
{\bf Proof:} Using the Cauchy-Schwartz inequality, we have:
\[
\begin{split}
\left|\textbf{I}(W_m(t))(z)-
\textbf{I}(W_m(s))(z)\right|&=\left|\sum_{k=1}^\infty
 \left\{(T_m\widetilde{h}_k)(t)-
 (T_m\widetilde{h}_k)(s)\right\}z_k\right|\\
&=\left|\sum_{k=1}^\infty \left\{(T_m\widetilde{h}_k)(t)-
(T_m\widetilde{h}_k)(s)\right\}
(2k)^{-\frac{q}{2}}(2k)^{\frac{q}{2}}z_k\right|\\
&\leq \left(\sum_{k=1}^\infty
\left|\left\{(T_m\widetilde{h}_k)(t)-
(T_m\widetilde{h}_k)(s)\right\}
\right|^2(2k)^{-q}\right)^\frac{1}{2}\times\\
&\hspace{5mm}\times \left(\sum_{k=1}^\infty \left|z^{\epsilon_k}
\right|^2 (2k)^q\right)^\frac{1}{2},
\end{split}
\]
and thus
\[
\begin{split}
\left|\textbf{I}(W_m(t))(z)- \textbf{I}(W_m(s))(z)\right|& \leq
\left(\sum_{k=1}^\infty \left|\left\{ (T_m\widetilde{h}_k)(t)-
(T_m\widetilde{h}_k)(s)\right\}
\right|^2(2k)^{-q}\right)^{\frac{1}{2}}\times\\
&\hspace{5mm}\times \left(\sum_{\alpha\in\ell}
\left|z^\alpha\right|^2
(2\N)^{\alpha q}\right)^\frac{1}{2}\\
&\leq \left|t-s\right|\left(\sum_{k=1}^\infty
 \left\{
\widetilde{C}_1k^\frac{11}{12} +\widetilde{C}_2
 \right\}^2(2k)^{-q}\right)^\frac{1}{2}\times\\
 &\hspace{5mm}\times
 \left(\sum_{\alpha\in\ell} \left|z\right|^{\alpha}
 (2\N)^{\alpha q}\right)^\frac{1}{2},
\end{split}
\]
where we have used \eqref{eq:uni_conv} to go from from the first
inequality to the second.\\

The first sum converges for $q\ge4$ and the second
converges since $z\in K_4(\delta)$, so we
conclude that $\textbf{I}(W_m(t))(z)$ is continuous
in $t$ for every $z\in K_4(\delta)$.
\mbox{}\qed\mbox{}\\

We now recall the following result of
\cite{MR1408433}, called {\sl the differentiation
of $S_{-1}$ processes}.

\begin{Pn}
\label{diff}
\cite[Lemma 2.8.4 p.
77]{MR1408433} Suppose $\left\{X(t,\omega)\right\}$ and
$\left\{F(t,\omega)\right\}$ are $S_{-1}$-valued processes such
that
\[
\frac{{ d}({\mathbf I}(X)(t))(z)}{{dt}}=({\mathbf I}(F)(t))(z)
\]
for each $t\in (a,b),~z\in K_q(\delta)$ and that $({\mathbf I}(F)(t))(z)$ is
bounded for $(t,z)\in (a,b)\times K_q(\delta)$, and is a
continuous function of $t$ for every $z\in K_q(\delta)$. Then
$X(t,\omega)$ is a differentiable process and
\[
\frac{dX(t,\omega)}{dt}=F(t,\omega)
\]
for all $t\in (a,b)$.
\end{Pn}

In view of this proposition, the first step toward showing that
$W_m$ is the derivative of $B_m$ is to show that this fact holds for the Hermite
transforms. This is done in the following lemma.

\begin{La}
\label{la:7.6}
Let $t\in\R$ and $z\in K_4(\delta)$. Then
$$\frac{d\textbf{I}(B_m(t))(z)}{dt}
=\textbf{I}(W_m(t))(z).$$
\end{La}
{\bf Proof:}
Let $h\in\R$. Then
\[
\begin{split}
&\left|\frac{\textbf{I}
(B_m(t+h))(z)-\textbf{I}(B_m(t))(z)}{h}-
\textbf{I}(W_m(t))(z)\right|\\
&=\frac{1}{\left|h\right|}\left|
\sum_{k=1}^\infty \int_t^{t+h}
\left((T_m\widetilde{h}_k)(s)-
(T_m\widetilde{h}_k)(t)\right)ds z_k\right|\\
&=\frac{1}{\left|h\right|}\left|
\sum_{k=1}^\infty \int_t^{t+h}
\left((T_m\widetilde{h}_k)(s)-
(T_m\widetilde{h}_k)(t)\right)ds
(2k)^{-\frac{q}{2}}(2k)^{\frac{q}{2}} z_k\right|\\
&\leq \frac{1}{\left|h\right|}
\left(\sum_{k=1}^\infty
\left|\int_t^{t+h}
\left((T_m\widetilde{h}_k)(s)-
(T_m\widetilde{h}_k)(t)
\right)ds\right|^2(2k)^{-q}\right)^\frac{1}{2}\cdot\\
&\hspace{10mm}
\left(\sum_{k=1}^\infty(2k)^q \left|z^{\epsilon_k}
\right|^2 \right)^\frac{1}{2}\\
&\leq \frac{1}{\left|h\right|}\left(\sum_{k=1}^\infty
\int_t^{t+h}\left|(T_m\widetilde{h}_k)(s)-
(T_m\widetilde{h}_k)(t)\right|^2ds(2k)^{-q}
\right)^\frac{1}{2}\cdot\\
&\hspace{10mm} \left(\sum_{\alpha\in\ell} (2\N)^{q\alpha}
\left|z^{\alpha}\right|^2 \right)^\frac{1}{2},
\end{split}
\]
and therefore
\[
\begin{split}
&\left|\frac{\textbf{I}
(B_m(t+h))(z)-\textbf{I}(B_m(t))(z)}{h}-
\textbf{I}(W_m(t))(z)\right|\\
&\leq \frac{1}{\left|h\right|} \left(\sum_{k=1}^\infty
\int_t^{t+h}\left|t-s\right|^2ds\left\{\widetilde{C}_1k^\frac{11}{12}
+\widetilde{C}_2
\right\}^2(2k)^{-q}\right)^\frac{1}{2}\cdot\\
&\hspace{10mm} \left(\sum_{\alpha\in\ell} (2\N)^{q\alpha}
\left|z^{\alpha}\right|^2 \right)^\frac{1}{2}\\
&\leq \frac{\left|h\right|^\frac{3}{2}}{\sqrt{3}
\left|h\right|}\left(\sum_{k=1}^\infty \left\{
\widetilde{C}_1k^\frac{11}{12} +\widetilde{C}_2 \right\}^2(2k)^{-q}
\right)^\frac{1}{2}\cdot\\
&\hspace{10mm} \left(\sum_{\alpha\in\ell} (2\N)^{q\alpha}
\left|z^{\alpha}\right|^2 \right)^\frac{1}{2}
\longrightarrow_{\left|h\right|\to 0}0.
\end{split}
\]
\mbox{}\qed\mbox{}

We are now ready for the main result of this
section:

\begin{Tm}
It holds that \[
 \frac{dB_m(t)}{dt}=W_m(t)
 \]
in the sense that
\[
 \frac{
d({\mathbf I}(B_m(t))(z))}{dt}={\mathbf I}(W_m(t))(z)
\]
for all $t\in\R$, point-wise boundedly.\\
\end{Tm}

{\bf Proof:} Taking in Proposition \ref{diff}, $X(t,\omega)=
B_m(t),~F(t,\omega)=W_m(t)$ we have from Lemma \ref{la:7.6}
\[
\frac{d(\textbf{I}(B_m(t))(z))}{dt}=\textbf{I}(W_m(t)(z))
\]
for all $(t,z)\in\R\times K_4(\delta)$, and by Proposition
\ref{pn:bd}, $\textbf{I}(W_m(t))(z)$ is a bounded function for all
$(t,z)\in\R\times K_2(\delta)$, then for all $(t,z)\in\R\times
K_4(\delta)$ the pair $(\textbf{I}(B_m(t))(z),\textbf{I}(W_m(t)))(z)$
satisfied the condition of Proposition \ref{diff} then we can
conclude that $B_m(t)$ is a differentiable in $S_{-1}$ process
which completes the proof.
\mbox{}\qed\mbox{}\\

\bibliographystyle{plain}
\def\cprime{$'$} \def\lfhook#1{\setbox0=\hbox{#1}{\ooalign{\hidewidth
  \lower1.5ex\hbox{'}\hidewidth\crcr\unhbox0}}} \def\cprime{$'$}
  \def\cprime{$'$} \def\cprime{$'$} \def\cprime{$'$} \def\cprime{$'$}

\end{document}